\documentclass[11pt]{amsart}

\usepackage[utf8]{inputenc}
\usepackage[T1]{fontenc}
\usepackage{amsmath, amssymb, amsthm}
\usepackage{mathrsfs}
\usepackage{enumitem}
\usepackage{hyperref}
\usepackage{geometry}
\usepackage{amsmath, amssymb}
\usepackage{algorithm}
\usepackage{amsthm}
\usepackage{amsmath, amssymb}
\usepackage{algorithm}
\usepackage[noend]{algpseudocode}
\usepackage{listings}
\usepackage{xcolor}
\usepackage{listings}
\usepackage{xcolor}
\usepackage{tcolorbox}
\usepackage{listings}
\usepackage{tcolorbox}
\usepackage{listings}
\newcommand{\multiset}[2]{\left( \genfrac{}{}{0pt}{}{#1}{#2} \right)}
\usepackage{xcolor}

\definecolor{lightgray}{gray}{0.95}
\definecolor{mathematicaBlue}{RGB}{0, 0, 128}

\lstdefinelanguage{Mathematica}{
  morekeywords={
    Module, Simplify, ClearAll, Print, Return, Table, TrueQ, D,
    If, And, Or, Not, True, False
  },
  sensitive=true,
  morecomment=[l]{(*},
  morestring=[b]"
}

\tcbset{colframe=black!60, colback=gray!10, fonttitle=\bfseries}

\lstdefinestyle{mmaoutput}{
  basicstyle=\ttfamily\small,
  backgroundcolor=\color{gray!10},
  frame=single,
  framerule=0pt,
  rulecolor=\color{gray!70},
  breaklines=true,
  showstringspaces=false,
  columns=flexible,
  literate={->}{{$\rightarrow$}}1
           {<=}{{$\le$}}1
           {>=}{{$\ge$}}1
           {==}{{$=$}}1
           {!=}{{$\neq$}}1
           {^}{{$^{\wedge}$}}1
           {_}{{$\_$}}1
}

\lstdefinelanguage{Mathematica}{
  morekeywords={
    Module, Simplify, ClearAll, Print, Return, Table, TrueQ, D,
    If, And, Or, Not, True, False
  },
  sensitive=true,
  morecomment=[l]{(*},
  morecomment=[s]{(*}{*)},
  morestring=[b]",
}

\newtheorem{theorem}{Theorem}[section]

\theoremstyle{definition}

\usepackage{listings}
\usepackage{xcolor}

\lstdefinestyle{matlab}{
  language=Matlab,
  basicstyle=\ttfamily\small,
  keywordstyle=\color{blue},
  commentstyle=\color{gray},
  stringstyle=\color{red},
  breaklines=true,
  frame=single,
  backgroundcolor=\color{gray!10},
  captionpos=b
}

\lstdefinelanguage{Mathematica}{
  morekeywords={
    Module, If, Break, For, True, False,
    Min, Length, Missing
  },
  sensitive=true,
  morecomment=[l]{(*},
  morestring=[b]",
}
\usepackage{pifont}

\lstdefinestyle{mathematica}{
  language=Mathematica,
  basicstyle=\ttfamily\small,
  keywordstyle=\color{blue},
  commentstyle=\color{gray},
  stringstyle=\color{red},
  breaklines=true,
  frame=single,
  backgroundcolor=\color{gray!10},
  captionpos=b
}

\usepackage{xcolor}

\theoremstyle{remark}

\title[On recurrence coefficients of classical orthogonal polynomials]{On recurrence coefficients of classical orthogonal polynomials}

\author{K. Castillo}
\address{CMUC, Department of Mathematics, University of Coimbra, 3000-143 Coimbra,
Portugal}
\email{kenier@mat.uc.pt}

\author{G. Gordillo-N\'u\~nez}
\address{CMUC, Department of Mathematics, University of Coimbra, 3000-143 Coimbra,
Portugal}
\email{up202310693@up.pt}

\dedicatory{}
\date{\today}

\begin{document}

\begin{abstract}
In \textit{Lett. Math. Phys.} \textbf{114}, 54 (2024) and \textbf{115}, 70 (2025), the author introduces what is presented as a novel method for determining whether a sequence of orthogonal polynomials is ``classical'', based solely on its initial recurrence coefficients. This note demonstrates that all the results contained in those works are already encompassed by two general theorems previously established in \textit{J. Math. Anal. Appl.} \textbf{515} (2022), Article 126390. A symbolic algorithm, implemented in \textsc{Mathematica}, is also provided to enable automated verification of the classical character of orthogonal polynomial sequences on quadratic lattices. As an application, it is shown that the so-called para-Krawtchouk polynomials on bi-lattices, discussed in \textit{Lett. Math. Phys.} \textbf{115}, 70 (2025), constitute a particular instance of a classical orthogonal family on a linear lattice. Consequently, their algebraic properties follow as a specific case of one of the main theorems established in \textit{J. Math. Anal. Appl.} \textbf{515} (2022), Article 126390.
\end{abstract}
\maketitle

\section{Introduction}
To ensure self-containment, the essential elements required to follow the core developments are briefly reviewed below, with topological considerations and other technical
subtleties intentionally omitted. A quick reading of Chapters 1 and 9 of \cite{CP25} ---which lay the foundations of the so-called algebraic theory of orthogonal polynomials, originally developed by Maroni ---may offer the reader a deeper perspective on these topics. Nonetheless, such familiarity is not strictly necessary, given the specific focus of the present work.

A sequence of polynomials \((P_n)_{n\in\mathbb{N}}\), with each \(P_n\) having exact degree \(n\), is said to be
orthogonal with respect to a linear functional \(\mathbf{u}\) on \(\mathbb{C}[\cdot]\), the space of all polynomial functions (or simply polynomials---in our context, we shall make no distinction between the two) with complex coefficients, whenever
\begin{align}\label{orth}
\langle \mathbf{u}, P_n\, P_m \rangle 
\begin{cases}
\neq 0, & \text{if } n = m, \\[7pt]
= 0, & \text{if } n \neq m.
\end{cases}
\end{align}
In a completely analogous fashion, one may construct a finite sequence of orthogonal polynomials satisfying the orthogonality relations for all degrees \( n \in \mathbb{N} \) such that \( n \leq M \), for some fixed natural number \( M \). A linear functional \( \mathbf{v}: \mathbb{C}[\cdot] \to \mathbb{C} \) is said to be {regular} if there exists a sequence of orthogonal polynomials with respect to \( \mathbf{v} \). In this context, a lattice is defined as any mapping \( x: \mathbb{C} \to \mathbb{C} \) of the form \cite[(1.9)]{CP25}
\begin{align}\label{lattices}
x(s) = 
\begin{cases}
\mathfrak{c}_1 q^{-s} + \mathfrak{c}_2 q^{s} + \mathfrak{c}_3, &q \neq 1, \\[7pt]
\mathfrak{c}_4 s^2 + \mathfrak{c}_5 s + \mathfrak{c}_6, & q = 1,
\end{cases}
\end{align}
where \( q \) is a fixed positive parameter and \( \mathfrak{c}_1, \dots, \mathfrak{c}_6 \) are complex constants such that \( (\mathfrak{c}_1, \mathfrak{c}_2) \neq (0,0) \) for \( q \neq 1 \) and \( (\mathfrak{c}_4, \mathfrak{c}_5, \mathfrak{c}_6) \neq (0,0,0) \) for \( q = 1 \). Observe that
\[
\displaystyle x\left(s + \frac{1}{2}\right) + x\left(s - \frac{1}{2}\right) = 2\,\alpha\, x(s) + 2\,\beta,
\]
where \(\alpha\) and \(\beta\) are given by
\begin{align*}
\alpha = \frac{q^{1/2} + q^{-1/2}}{2}, \quad
\beta = \begin{cases}
(1 - \alpha)\, \mathfrak{c}_3, & q \neq 1, \\[7pt]
\displaystyle \frac{1}{4}\, \mathfrak{c}_4, & q = 1.
\end{cases}
\end{align*}

Given a lattice \( x \) and a polynomial \( p \), we define the $x$-derivative and the $x$-average \cite[Definition 1.2]{CP25}, denoted respectively by \( \mathrm{D} \) and \( \mathrm{S} \), as the linear operators satisfying $\deg(\mathrm{D} p)=\deg p-1$, $\deg(\mathrm{S} p)=\deg p$, and
\begin{align*}
\mathrm{D} p(x(s))&= \frac{p\left(x\left(s+\displaystyle \frac12\right)\right)-p\left(x\left(s-\displaystyle \frac12\right)\right)}{x\left(s+\displaystyle \frac12\right)-x\left(s-\displaystyle \frac12\right)}, \\[7pt]
\mathrm{S} p(x(s))&= \frac{1}{2}\,p\left(x\left(s+\displaystyle \frac12\right)\right)-\frac{1}{2}\,p\left(x\left(s-\displaystyle \frac12\right)\right),
\end{align*}
where $\mathrm{D} p=p'$ and  $\mathrm{S} p=p$ whenever $x(s)=\mathfrak{c}_6$. The reader should bear in mind that
these operators depend intrinsically on the fixed lattice, even though their dependence on
the ``variable'' is often suppressed for simplicity---much like it is occasionally done in the
context of the standard derivative. By defining, via transposition, certain operations on
a given linear functional \(  \mathbf{v}: \mathbb{C}[\cdot] \to \mathbb{C} \)---such as multiplication by a polynomial \(p\) \cite[Definition 1.1]{CP25}, \(x\)-differentiation, and \(x\)-averaging \cite[Definition 1.3]{CP25}---
denoted respectively by \(p\mathbf{v}\), \(\mathbf{Dv}\), and \(\mathbf{Sv}\), and by applying standard principles
of duality, it can be shown that all the families in the Askey scheme of hypergeometric
orthogonal polynomials and their q-analogues are orthogonal with respect to a regular
linear functional \(\mathbf{u}: \mathbb{C}[\cdot] \to \mathbb{C} \) satisfying the equation
\begin{align}\label{pearson}
\mathbf{D}(\phi\, \mathbf{u}) = \mathbf{S}(\psi\, \mathbf{u}),
\end{align}
where \( \phi \) is a nonzero polynomial of degree at most two, and \( \psi \) is a polynomial of degree exactly one, both independent of \( n \). In this note, we follow the definition given in \cite[Definition 3.1]{CMP22a} and \cite[Definition 9.1]{CP25}, according to which a sequence of orthogonal polynomials with respect to  \(\mathbf{u}\)  is said to be classical if it satisfies the equation \eqref{pearson}. A detailed understanding of the origin of \eqref{pearson} is not required to follow the present note; the reader may simply rely on the previously stated assertions and the theorems that follow. 

\begin{theorem}\label{thm1}\cite[Theorem 4.1]{CMP22a}\cite[Theorem 9.2]{CP25}
Consider the lattice
\[
x(s) = \mathfrak{c}_1 q^{-s} + \mathfrak{c}_2 q^{s} + \mathfrak{c}_3.
\]
Let \({\bf u} : \mathbb{C}[\cdot] \to \mathbb{C}\) be a nontrivial linear functional, and suppose that there exist polynomials \(\phi\) and \(\psi\), of degree at most two and one respectively, not both identically zero, such that
\[
\mathbf{D}(\phi\,{\bf u}) = \mathbf{S}(\psi\,{\bf u}).
\]
Let
\[
\phi(z) = a z^2 + b z + c, \quad \psi(z) = d z + e.
\]
Then, \({\bf u}\) is regular if and only if
\begin{equation}
\label{regq}
    d_n \neq 0 \quad \text{and} \quad \phi^{[n]}\left(\mathfrak{c}_3 - \frac{e_n}{d_{2n}}\right) \neq 0, \quad \text{for all } n \in \mathbb{N},
\end{equation}
where
\begin{align*}
\phi^{[n]}(z)&=\left(d(\alpha^2-1)\gamma_{2n}+a\alpha_{2n}\right)\left((z-\mathfrak{c}_3)^2-2\mathfrak{c}_1\mathfrak{c}_2\right)  \\[7pt]
&\quad +\left(\phi'(\mathfrak{c}_3)\alpha_n+\psi(\mathfrak{c}_3)(\alpha^2-1)\gamma_n\right)(z-\mathfrak{c}_3)
+ \phi(\mathfrak{c}_3)+2a\mathfrak{c}_1\mathfrak{c}_2,
\end{align*}
and
\begin{align*}
d_n= a \gamma_n + d \alpha_n,\quad e_n=(2a \mathfrak{c}_3 + b) \gamma_n + (d \mathfrak{c}_3 + e) \alpha_n,
\end{align*}
with
\[
\alpha_n = \frac{q^{n/2} +q^{-n/2}}{2},\quad \gamma_n=\displaystyle\frac{q^{n/2}-q^{-n/2}}{q^{1/2}-q^{-1/2}}.
\]
Only under the aforementioned regularity conditions does there exist a monic orthogonal polynomial sequence \((P_n)_{n \in \mathbb{N}}\)---possibly only up to some finite degree, if condition~\eqref{regq} is satisfied for \( n \in \mathbb{N} \) with \( n \leq M \), for some fixed natural number \( M \)---that satisfies the three-term recurrence relation
\[
P_{n+1}(z) = (z - B_n)\, P_n(z) - C_n\, P_{n-1}(z),
\]
with initial conditions \(P_{-1}(z) = 0\), \(P_0(z) = 1\), and recurrence coefficients given by
\begin{align*}
B_n & =\mathfrak{c}_3+ \frac{\gamma_n e_{n-1}}{d_{2n-2}}-\frac{\gamma_{n+1}e_n}{d_{2n}}, \\[7pt]
C_{n+1} & =-\frac{\gamma_{n+1}d_{n-1}}{d_{2n-1}d_{2n+1}}\phi^{[n]}\left(\mathfrak{c}_3 -\frac{e_{n}}{d_{2n}}\right).
\end{align*}
\end{theorem}

\begin{theorem}\label{thm2}\cite[Theorem 4.2]{CMP22a}\cite[Theorem 9.3]{CP25}
Consider the quadratic  lattice
\[
x(s) = \mathfrak{c}_4 s^2 + \mathfrak{c}_5 s + \mathfrak{c}_6.
\]
Let \({\bf u} : \mathbb{C}[\cdot] \to \mathbb{C}\) be a nontrivial linear functional, and suppose that there exist polynomials \(\phi\) and \(\psi\), of degree at most two and one respectively, not both identically zero, such that
\[
\mathbf{D}(\phi\,{\bf u}) = \mathbf{S}(\psi\,{\bf u}).
\]
Define
\[
\phi(z) = a z^2 + b z + c, \quad \psi(z) = d z + e.
\]
Then, \({\bf u}\) is regular if and only if
\begin{equation}
\label{regqua}
    d_n \neq 0, \quad \phi^{[n]}\left(-\frac{1}{4} \mathfrak{c}_4 n^2 - \frac{e_n}{d_{2n}} \right) \neq 0, \quad \text{for all } n \in \mathbb{N},
\end{equation}
where
\[
d_n = a n + d, \quad
e_n = b n + e + \frac{1}{2} \mathfrak{c}_4 d n^2,
\]
and
\begin{align*}
\phi^{[n]}(z) &= a z^2 + \left(b + \frac{3}{2} \mathfrak{c}_4 n d_n \right) z + \phi\left(\frac{1}{4} \mathfrak{c}_4 n^2 \right) + \frac{1}{2} \mathfrak{c}_4 n \psi\left(\frac{1}{4} \mathfrak{c}_4 n^2 \right)\\[7pt]
&\quad - \frac{n}{4} \left(16 \mathfrak{c}_4 \mathfrak{c}_6 - \mathfrak{c}_5^2\right) d_n.
\end{align*}
Only under the aforementioned regularity conditions does there exist a monic orthogonal polynomial sequence \((P_n)_{n \in \mathbb{N}}\)---possibly only up to some finite degree, if condition \eqref{regq} is satisfied for \( n \in \mathbb{N} \) with \( n \leq M \), for some fixed natural number \( M \)---that satisfies the three-term recurrence relation
\[
P_{n+1}(z) = (z - B_n)\, P_n(z) - C_n\, P_{n-1}(z),
\]
with initial conditions \(P_{-1}(z) = 0\), \(P_0(z) = 1\), and recurrence coefficients given by
\begin{align*}
B_n &= \frac{n e_{n-1}}{d_{2n-2}} - \frac{(n+1) e_n}{d_{2n}} - \mathfrak{c}_4 n(n - 1), \\[6pt]
C_{n+1} &= -\frac{(n+1) d_{n-1}}{d_{2n-1} d_{2n+1}}\, \phi^{[n]}\left(-\frac{1}{4} \mathfrak{c}_4 n^2 - \frac{e_n}{d_{2n}} \right).
\end{align*}\end{theorem}

In Section \ref{2}, it will be shown, in a single step, that all the results of \cite{M24, M25} are already contained
in Theorems \ref{thm1} and \ref{thm2}. Section \ref{3} presents an algorithm for determining whether a given pair of sequences, interpreted as the recurrence coefficients of a monic orthogonal polynomial sequence, corresponds to a classical family on  \(x(s) = \mathfrak{c}_4 s^2 + \mathfrak{c}_5 s + \mathfrak{c}_6\). Section \ref{4}
demonstrates ---both formally and through the algorithm presented in the preceding section ---
that the so-called para-Krawtchouk polynomials on bi-lattices are simply a special case of classical orthogonal polynomials on a linear lattice, though conveyed in an inadvertently obscure and intricate manner. Consequently, their regularity, recurrence coefficients, and all associated algebraic properties follow directly from previously established theorems, rendering the main theorem of \cite{CFM25} largely superfluous.

\section{A ``Proof'' of the Main Results of \cite{M24, M25}}\label{2}

Assume that \((P_n)_{n \geq 0}\) is a sequence of classical orthogonal polynomials  in the sense of \eqref{orth}. Consequently, the recurrence coefficients are determined by Theorem \ref{thm1} or Theorem \ref{thm2}. By evaluating \(B_n\) and \(C_{n+1}\) at \(n = 0\) and \(n = 1\) in these theorems, one obtains explicit expressions for the polynomials \(\phi\) and \(\psi\) in terms of \(B_0\), \(B_1\), \(C_1\), and \(C_2\). (It is worth noting that the dependence of \( \phi \) and \( \psi \) on the recurrence coefficients \( B_0 \), \( B_1 \), \( C_1 \), and \( C_2 \) is well known to those familiar with this line of work initiated by Maroni (as surveyed in \cite{CP25}); in any case, this dependency is explicitly addressed in \cite[Lemma 2.2]{CMP22b}.) This elementary observation underlies the analysis carried out in \cite{M24,M25} ---despite the fact that the author offers long and circuitous proofs, starting from a definition of classical orthogonal polynomials that differs from that of \cite{CMP22a, CP25}, only to ultimately depend on the same definition and the essential role played by the polynomials \(\phi\) and \(\psi\). While \cite{CMP22b, CMP22a} are cited in \cite{M24, M25}, their contributions are not acknowledged in the key conceptual passages where they would be most relevant. As a result, the reader may be left with the mistaken impression that certain results or techniques are novel, when in fact they had already been established in earlier works and even included in recent monographs on the subject.

\subsection{Case \(x(s) = \mathfrak{c}_4 s^2 + \mathfrak{c}_5 s + \mathfrak{c}_6\)} Let us apply Theorem \ref{thm2} and assume, without loss of generality, that \(d = 1\). If, on the contrary, \(d = 0\), then, by Theorem \ref{thm2}, the functional \(\mathbf{u}\) would not be regular. From the expression for \(B_n\), it follows immediately that \(B_0 = -e\), simply by evaluating at \(n = 0\). This suffices to fully determine the polynomial \(\psi\).  Next, we compute \(B_1\), \(C_1\), and \(C_2\) by evaluating \(B_n\) at \(n = 1\) and \(C_{n+1}\) at \(n = 0\) and \(n = 1\). This yields an elementary system of equations in the unknown coefficients \(a\), \(b\), and \(c\) of the polynomial \(\phi\), whose solution yields the following explicit expressions::
\[
a = \mathfrak{a}, \quad b = -\mathfrak{a} B_0 - \mathfrak{b}, \quad c = \mathfrak{b} B_0 - (\mathfrak{a} + 1) C_1,
\]
where
\begin{align*}
\mathfrak{a} &= -\frac{1}{3} + \frac{(B_1 - B_0)^2 - 8\beta(B_0 + B_1 - 2\beta) + 4(C_1 + 4\beta c_6 - c_5^2)}{6C_2},\\[7pt]
\mathfrak{b} &= (\mathfrak{a} + 1) B_1 + 2\beta - \frac{B_0 + B_1}{2}.
\end{align*}
This is precisely the content of \cite[ Lemma 2.3]{M25}. As for \cite[ Theorem 2.4]{M25}, its statement becomes redundant in light of Theorem \ref{thm2}, which already provides explicit expressions for \(B_n\) and \(C_{n+1}\). Finally, \cite[ Theorem 2.7]{M25} is merely a reformulation of Theorem \ref{thm2} in a different guise. The reader may also compare \cite[ Lemma 2.3]{M25} and \cite[ Theorem 2.4]{M25} with their counterparts in \cite[ Lemma 2.2]{CMP22b} and \cite[Lemma 2.3]{CMP22b}, respectively, although the analysis in \cite{CMP22b} is more delicate, as it deals with a specific structural relation.

\subsection{Case \(x(s) = \mathfrak{c}_1 q^{-s} + \mathfrak{c}_2 q^{s} + \mathfrak{c}_3\)}
The particular case \( c_1 = c_2 =1/2 \) and \( c_3 = 0 \) was precisely the one considered in \cite{M24}. We leave it to the reader to verify, using Theorem~\ref{thm1}, both the results presented in \cite{M24} and their straightforward extensions.

\section{Identifying Classical Orthogonal Polynomials on \(x(s) = \mathfrak{c}_4 s^2 + \mathfrak{c}_5 s + \mathfrak{c}_6\)}\label{3}
To address aspects not covered in \cite{M25}, we present a symbolic algorithm, implemented in \textsc{Mathematica}, which will be used in the following section to identify the so-called para-Krawtchouk polynomials discussed in that work. The algorithm is designed to determine whether a given pair of sequences interpreted as the recurrence coefficients of a monic orthogonal polynomial sequence, corresponds to a classical family defined on the lattice
\(
x(s) = \mathfrak{c}_4 s^2 + \mathfrak{c}_5 s + \mathfrak{c}_6.
\) The algorithm, based entirely on Theorem \ref{thm2}, checks whether the sequences are the recurrence coefficients of a classical orthogonal polynomial family. If so, it returns the coefficients of the polynomials  
\[
\phi(z) = a z^2 + b z + c, \quad \psi(z) = z + e.
\]  
(Recall from the previous section that \(d\) may be taken to be equal to \(1\) without loss of generality.)
This is particularly useful in view of the following section, where we demonstrate that the para-Krawtchouk polynomials on bi-lattices discussed in \cite{M25} constitute a classical orthogonal family on the lattice, for instance\footnote{It is important to recall that the parameter \(\mathfrak{c}_6\) plays no essential role for linear lattices \cite{CMP22a}.}, \(x(s) = 2s + 1\). The reader may, when a genuine need arises, extend these ideas to the lattices of the form \(x(s) = \mathfrak{c}_1 q^{-s} + \mathfrak{c}_2 q^{s} + \mathfrak{c}_3\) and will find in Theorem~\ref{thm1} all the necessary ingredients to do so.\\
\begin{lstlisting}[language=Mathematica]
ClearAll["Global`*"];

formatXS[c4_, c5_, c6_] := Module[{terms = {}, expr},
  If[c4 =!= 0, AppendTo[terms, If[c4 === 1, "s^2", If[c4 === -1, "-s^2", ToString[c4] <> " s^2"]]]];
  If[c5 =!= 0, AppendTo[terms, If[c5 === 1, "s", If[c5 === -1, "-s", ToString[c5] <> " s"]]]];
  If[c6 =!= 0, AppendTo[terms, ToString[c6]]];
  expr = If[terms === {}, "0", StringRiffle[terms, " + "]];
  Return[expr];
];

ClassicalQuadratic[BN_, CNp1_, c4_, c5_, c6_] := 
 Module[{a, b, c, e, beta, d = 1, B0, B1, C1, C2, phi, psi, dphi, dN,
   eN, phiN, BNnew, CNp1new, someBNChecks, someCNp1Checks, 
   compareBN, compareCNp1, equals = False, xs},

  B0 = BN[0]; B1 = BN[1]; C1 = CNp1[0]; C2 = CNp1[1];
  beta = c4/4; e = -B0;

  a = Simplify[(B0^2 + 4 C1 - 2 C2 + (B1 - c4)^2 - 
       2 B0 (B1 + c4) - c5^2 + 4 c4 c6)/(6 C2)];

  b = Simplify[-((B0^3 + B1^3 - 2 B1^2 c4 + 3 C2 c4 - 
       B0^2 (B1 + 2 c4) + B1 (4 C1 + C2 + c4^2 - c5^2 + 4 c4 c6) - 
       B0 (B1^2 - 4 C1 + 5 C2 + 4 B1 c4 + c5^2 - 
          c4 (c4 + 4 c6)))/(6 C2))];

  c = Simplify[
    1/(6 C2) (B0^3 B1 - B0^2 (C1 + 3 C2 + 2 B1 (B1 + c4)) - 
      C1 (4 (C1 + C2) + (B1 - c4)^2 - c5^2 + 4 c4 c6) + 
      B0 (B1^3 - 2 B1^2 c4 + (2 C1 + 3 C2) c4 + 
         B1 (6 C1 + C2 + c4^2 - c5^2 + 4 c4 c6)))];

  phi[x_] := a x^2 + b x + c;
  psi[x_] := x + e;
  dphi[x_] := D[phi[x], x];
  dN[n_] := a n + d;
  eN[n_] := b n + e + 2 beta d n^2;
  phiN[x_, n_] := 
   a x^2 + (b + 6 beta n dN[n]) x + phi[beta n^2] + 
    2 beta n psi[beta n^2] - n/4 (16 beta c6 - c5^2) dN[n];

  BNnew[n_] := 
   Simplify[(n eN[n - 1])/dN[2 n - 2] - ((n + 1) eN[n])/dN[2 n] - 
     2 beta n (n - 1)];

  CNp1new[n_] := 
   Simplify[-((n + 1) dN[n - 1])/(dN[2 n - 1] dN[2 n + 1])*
     phiN[-beta n^2 - eN[n]/dN[2 n], n]];

  someBNChecks = Table[TrueQ[Simplify[BN[k] == BNnew[k]]], {k, 2, 5}];
  If[AllTrue[someBNChecks, TrueQ],
   someCNp1Checks = 
    Table[TrueQ[Simplify[CNp1[k] == CNp1new[k]]], {k, 2, 5}];
   If[AllTrue[someCNp1Checks, TrueQ],
    compareBN = TrueQ[Simplify[BN[n] == BNnew[n]]];
    If[compareBN,
     compareCNp1 = TrueQ[Simplify[CNp1[n] == CNp1new[n]]];
     If[compareCNp1, equals = True];
     ]
    ]
   ];

  xs = formatXS[c4, c5, c6];

  If[equals,
   Print["The sequence is classical for x(s) = ", xs];
   Print["a = ", a];
   Print["b = ", b];
   Print["c = ", c];
   Print["e = ", e],
   Print["The sequence is not classical for x(s) = ", xs]
   ];

  Return[equals];
];
\end{lstlisting}

\section{Who are the para-Krawtchouk polynomials?}\label{4}
Consider the polynomials discussed in \cite[Remark~2.6]{M25} and originally introduced in \cite{VZ12}, whose recurrence coefficients are given in \cite[p.~7]{CFM25} and, for  \(n < N\), take the form
\begin{align*}
B_n &= \frac{N +  \gamma - 1}{2}, \\[7pt]
C_{n+1} &= -\frac{(n+1)(n - N)(2n + 1 - N -  \gamma)(2n + 1 - N +  \gamma)}{4(2n - N)(2n - N + 2)},
\end{align*}
where \(N\) is a fixed odd positive integer and \( \gamma \) is a positive number less than \(2\).  Observe that, under these conditions, \(C_{n+1}\) remains strictly positive up to $C_N$ when it vanishes, thereby placing the problem within the finite positive-definite framework. However, in the broader context of regular linear functionals, one may even, for instance, consider \( N \in \mathbb{C} \setminus \mathbb{N} \) and allow \( \gamma \) to take complex values, provided the condition
\[
\gamma^2 - (2n + 1 - N)^2 \ne 0,
\]
is satisfied. Under these assumptions, the associated polynomials continue to exist and remain orthogonal with respecto a regular linear functional. However, within the scope of this note, our only intention is to convey to the reader that when the values of \( C_{n+1} \) are nonzero, the corresponding orthogonality holds with respect to a regular linear functional; and that, whenever these values are positive, such a functional becomes positive-definite, with all the associated consequences.

In \cite[Remark 2.6]{M25}, the author asserts that \cite{CFM25} established the classical nature of this family in the sense considered in the present note. Indeed, \cite[Theorem 1.1]{CFM25}---whose proof extends over approximately thirteen journal pages---is presented as an analogue of Theorem \ref{thm2}. However, this result is, in essence, nothing other than Theorem \ref{thm2} itself, specialized to the lattice \(x(s) = 2s + 1\). Any remaining doubt may be dispelled by setting \(\mathfrak{c}_4 = 0\), \(\mathfrak{c}_5 = 2\), and \(\mathfrak{c}_6 = 1\) in Theorem~\ref{thm2} and computing the corresponding expressions explicitly. The reason is simple: the para-Krawtchouk polynomials form a classical family with respect to any lattice of the form  \(x(s) = 2s +  \mathfrak{c}_6\). The algorithm from Section \ref{3} may now be employed to verify this claim. In the following box, the recurrence coefficients of the para-Krawtchouk polynomials are introduced. The values \(\mathfrak{c}_4 = 0\), \(\mathfrak{c}_5 = 2\), and \(\mathfrak{c}_6 = 1\) are assumed, and the \textsc{Mathematica} function \texttt{ClassicalQuadratic} is applied.

\begin{tcolorbox}[title=\textcolor{mathematicaBlue}{}, colback=lightgray, colframe=white, boxrule=0pt]\begin{lstlisting}[language=Mathematica]
BN[n_] := 1/2 (-1 + N + \[Mu]);
CNp1[n_] := ((1 + n) (n - N) (-(-1 - 2 n + N)^2 + \[Mu]^2))/
            (4 (2 n - N) (2 + 2 n - N));

(* Test for x(s) = 2s *)
ClassicalQuadratic[BN, CNp1, 0, 2, 1]
\end{lstlisting}
\end{tcolorbox}

The output of the \textsc{Mathematica} algorithm is:

\begin{tcolorbox}[title=\textcolor{mathematicaBlue}{}, colback=lightgray, colframe=white, boxrule=0pt]The sequence is classical for $x(s) = 2s+1$

\bigskip

\[
\begin{aligned}
a &= \frac{1}{1 - N} \\
b &= \frac{-1 + N +  \gamma}{-1 + N} \\
c &= \frac{1}{2}(1 - N -  \gamma) \\
e &= \frac{1}{2}(1 - N -  \gamma)
\end{aligned}
\]
\end{tcolorbox}

The reader should bear in mind that the computations performed in \textsc{Mathematica} are symbolic in nature, and that we always assume we are starting from a sequence of orthogonal polynomials. Consequently, the domain of the parameters involved is well defined from the outset. And what happens when the same question is posed with respect to the lattice  \(x(s) = 3s\)?

\begin{tcolorbox}[title=\textcolor{mathematicaBlue}{}, colback=lightgray, colframe=white, boxrule=0pt]\begin{lstlisting}[language=Mathematica]
(* Test for x(s) = 3s *)
ClassicalQuadratic[BN, CNp1, 0, 3, 0]
\end{lstlisting}
\end{tcolorbox}

The output of the \textsc{Mathematica} algorithm in this case is:

\begin{tcolorbox}[title=\textcolor{mathematicaBlue}{}, colback=lightgray, colframe=white, boxrule=0pt]The sequence is not classical for $x(s) = 3s$
\end{tcolorbox}

In light of the preceding discussion, let us attempt to interpret the para-Krawtchouk polynomials  through the lens of Theorem \ref{thm2}. What is referred to in \cite{VZ12} as a bi-linear lattice is a mapping \( y: \mathbb{N}_N = \{0,1,\dots,N\} \to \mathbb{R} \), for a fixed odd integer \(N\), of the form
\begin{align*}
    y(s)=s+\frac{1}{2}(\gamma-1)(1-(-1)^s),
\end{align*}
where \( \gamma \) is a positive number less than \(2\). We introduce the linear functional \({\bf w} : \mathbb{C}[\cdot] \to \mathbb{C}\), defined by
\[
\langle \mathbf{w}, p \rangle = \sum_{s \in \mathbb{N}_N} p(y(s))\, \omega(y(s)),
\]
where the weight function \( \omega \) is given by
\[
\omega(y(s)) = \mathbf{1}_{2\mathbb{N}}(y(s))\, \omega_1(y(s)) + \big(1 - \mathbf{1}_{2\mathbb{N}}(y(s))\big)\, \omega_2\big(y(s) - \gamma\big),
\]
and \( \mathbf{1}_{2\mathbb{N}} \) denotes the characteristic function of the even natural numbers. The functions \( \omega_1 \) and \( \omega_2 \) are real-valued and defined on \( 2\mathbb{N} \) by
\begin{align*}
    \omega_1(m) &= \frac{2^{-N} \left(1 + \displaystyle\frac{\gamma}{2}\right)_J}{\left(\displaystyle \frac{1}{2} \right)_J}  
    \frac{(-J)_{\tfrac{m}{2}} \left(-\displaystyle \frac{\gamma}{2} - J\right)_{\tfrac{m}{2}}}{\left(\displaystyle \frac{m}{2}\right)! \left(1 - \displaystyle \frac{\gamma}{2} \right)_{\tfrac{m}{2}}}, \\[8pt]
    \omega_2(m) &= \frac{2^{-N} \left(1 - \displaystyle\frac{\gamma}{2} \right)_J}{\left(\displaystyle \frac{1}{2} \right)_J}  
    \frac{(-J)_{\tfrac{m}{2}} \left(\displaystyle \frac{\gamma}{2} - J \right)_{\tfrac{m}{2}}}{\left(\displaystyle \frac{m}{2}\right)! \left(1 + \displaystyle \frac{\gamma}{2} \right)_{\tfrac{m}{2}}},
\end{align*}
with \( 2J = N - 1 \). Observe that \(\mathbf{w}\) defined in this way is precisely the one with respect to which the para-Krawtchouk polynomials \((P^{(\gamma)}_n)_{n\in \mathbb{N}_N}\)  are orthogonal \cite[Section 3]{VZ12}. From this point onward, the only requirement for the reader is to recall that, in Theorem~\ref{thm2}, the lattice is, in general, a mapping from a subset of \( \mathbb{C} \) into \( \mathbb{C} \). One need only observe that, if we define a function \( x : \mathbb{V}_N \to \mathbb{R} \) by
\[
x(s)=2s+1,
\]
with
\[
\mathbb{V}_N = \left\{-\frac{1}{2},\ \frac{\gamma}{2} - \frac{1}{2},\ \frac{1}{2},\ \frac{\gamma}{2} + \frac{1}{2},\ \dots,\ \frac{N-2}{2},\ \frac{\gamma}{2} + \frac{N-2}{2} \right\},
\]
then
\[
\displaystyle
y(\mathbb{N}_N) = \{0,\gamma,2,\gamma+2,\dots,N-1,\gamma+N-1\}=x(\mathbb{V}_N).
\]
Hence, the orthogonality of \((P_n^{(\gamma)})_{n\in \mathbb{N}_N}\) with respect to the functional \( \mathbf{w} \) may be equivalently written as
\begin{align*}
\langle \mathbf{w}, P^{(\gamma)}_n\,P^{(\gamma)}_m \rangle&=\sum_{s\in \mathbb{N}_N} P^{(\gamma)}_n(y(s)) P^{(\gamma)}_m(y(s)) \omega(y(s))\\[7pt]
&=\sum_{s\in \mathbb{V}_N} P^{(\gamma)}_n(x(s)) P^{(\gamma)}_m(x(s)) \omega(x(s)).
\end{align*}
This merely serves to confirm what was already established by the algorithm presented in the previous section.

\section*{Conflict of Interest}

On behalf of all authors, the corresponding author states that there is no conflict of interest.

\section*{Data Availability Statement}

This article does not contain any datasets. All results and conclusions are based on theoretical analysis and mathematical reasoning.
\section*{Acknowledgements}
This work was partially supported by the Centre for Mathematics of the University of Coimbra, funded by the Portuguese Government through FCT/MCTES (DOI: 10.54499/ UIDB/00324/2020). The first author acknowledges financial support from FCT under the grant DOI: 10.54499/2022.00143.CEECIND/CP1714/CT0002. The second author acknowledges financial support from FCT under the grant DOI: 10.54499/UI.BD. 154694.2023.

\bibliographystyle{plain}      

\bibliography{bib}  

\end{document}